\documentclass[letterpaper,10 pt,conference]{ieeeconf} 

\IEEEoverridecommandlockouts

\usepackage{amsmath,amssymb,verbatim}
\usepackage{color,balance}
\usepackage{cgmMacros,mathfontmacros}

\renewcommand{\Re}{\mathbb{R}}
\newcommand{\SC}{\mathbb{S}}
\newcommand{\A}{\mathcal{A}}
\newcommand{\Y}{\mathcal{Y}}
\newcommand{\Kinf}{\mathcal{K}_{\infty}}
\newcommand{\X}{\mathbb{X}}

\newcommand{\D}{\mathcal{D}}
\newcommand{\E}{\mathcal{E}}
\newcommand{\W}{\mathcal{W}}
\newcommand{\B}{\mathcal{B}}

\newcommand{\tcb}{\textcolor{blue}}

\definecolor{dgreen}{rgb}{0,0.4,0}

\newtheorem{theorem}{Theorem}

\newtheorem{lemma}{Lemma}
\newtheorem{example}{Example}
\newtheorem{proposition}{Proposition}
\usepackage{setspace}
\begin{document}

\title{Further results on synergistic Lyapunov functions \\ and hybrid
  feedback design through backstepping \vspace{-0.0in}}

\author{Christopher G. Mayhew, Ricardo G. Sanfelice, and Andrew
  R. Teel \vspace{-0.0in}
\thanks{Dr. Mayhew is with the Robert Bosch Research and
    Technology Center, Palo Alto, CA, 94304. Dr. Sanfelice is with
    Department of Aerospace and Mechanical Engineering, University of
    Arizona, Tucson, AZ 85721-0119.  Dr. Teel is with the ECE
    Department, University of California, Santa Barbara, CA
    93106-9560.  Research supported in part by the AFOSR grant
    FA9550-09-1-0203 and the NSF grants ECCS-0925637, and
    CNS-0720842.}  }

\maketitle

\thispagestyle{empty}
\pagestyle{empty}

%
%


\begin{abstract}
We extend results on backstepping hybrid feedbacks by exploiting
synergistic Lyapunov function and feedback (SLFF) pairs in a
generalized form.  Compared to existing results, we delineate SLFF
pairs that are ``ready-made'' and do not require extra dynamic
variables for backstepping.  From an (weak) SLFF pair for an affine
control system, we construct an SLFF pair for an extended system where
the control input is produced through an integrator.  The
resulting hybrid feedback asymptotically stabilizes the extended
system when the ``synergy gap'' for the original system is strictly
positive.  To highlight the versatility of SLFF pairs, we provide a
result on the existence of a SLFF pair whenever a hybrid feedback
stabilizer exists. The results are illustrated 
on
the ``3D pendulum.''
\end{abstract}

\vspace{-.1in}

\section{Introduction}
Hybrid feedbacks are commonly used to improve performance and achieve
objectives that elude classical feedback designs.  Such objectives
include global asymptotic stabilization of a point for a system
evolving on a manifold that is not topologically equivalent to a
Euclidean space, or global asymptotic stabilization of a disconnected
set of points.

In a recent series of results, synergistic potential functions are
developed and used to achieve robust, global asymptotic stabilization
of planar orientation \cite{Mayhew2010-planar}, orientation on the
2-sphere \cite{Mayhew2010-spherical} (applied to the 3D pendulum in
\cite{Mayhew2010-3dpend}), and rigid-body attitude
\cite{Mayhew2011-so3potentials, Mayhew2011-so3control}.  Synergistic
potential functions are extended to synergistic Lyapunov function and
feedback (SLFF) pairs in \cite{Mayhew2011-synergistic}.  For a
continuous-time system with embedded logic variables, a Lyapunov
function and feedback pair is synergistic when, at places in the
state-space where the feedback is ineffective, the logic variable can
be switched to decrease the value of the Lyapunov function.  The
magnitude of the available decrease is called the {\em synergy gap}.
In \cite{Mayhew2011-synergistic}, the synergy gap is defined as an
infimum over an appropriate subset of the state space and it is
required to be positive for control synthesis and backstepping.  In
this note, the synergy gap is state dependent. It must be positive
away from a desired compact set and everywhere positive for
backstepping that achieves global asymptotic stability.

Earlier control algorithms propose a similar scheme exploiting
multiple Lyapunov functions. Some have appeared in the context of
adaptive control using hysteresis \cite{Morse1992,Hespanha1999} and
supervisory control systems \cite{Malmborg1996}.  Applications using
this feedback scheme have appeared for swing-up and stabilization of
an inverted pendulum \cite{Malmborg1996, Fierro1999} and for control
of a double-tank system \cite{Malmborg1997}. Multiple Lyapunov
functions are also proposed for control and analysis in
\cite{Branicky1998}.

In Section~\ref{section:HC} we construct a robustly globally
asymptotically stabilizing hybrid feedback algorithm using an SLFF
pair. In Sections \ref{section:BS0}--\ref{section:BS4}, we broaden the
applicability of (weak) SLFF pairs through backstepping.  Starting
from a weak SLFF pair for an affine control system, we construct an
(non-weak) SLFF pair for an extended system where the control is
produced through an integrator.   Results of this type for continuous-time systems 
can be found in
\cite[Lemma 2.8(ii)]{KKKBook} and
\cite[Theorem 5.3]{Byrnes1991a}.  Similar results for
switched systems appear in \cite{Xiang2008}; however, the
notion of synergism that is crucial for ensuring global asymptotic
stability does not appear in \cite{Xiang2008}.
%
%
We provide a variety of backstepping results:

$\bullet$ The backstepping algorithm resembles classical
backstepping when the assumed weak SLFF pair is
{\em pure} and is {\em ready-made} relative to a quadratic function.
An SLFF pair is {\em pure} when the Lyapunov function is
non-increasing along solutions at every point in the state space when
using the feedback.  An SLFF pair is {\em ready-made} when there is an
appropriate relationship between the size of the jumps in the feedback
law and the synergy gap of the SLFF pair.  These definitions are made
precise in Section \ref{section:SLFF2} and the backstepping algorithm
is described in Section \ref{section:BS1}.

$\bullet$ 
If the weak SLFF pair is not pure, a backstepping result can still be
obtained when the SLFF is ready-made relative to a linear function.
See the algorithm in Section \ref{section:BS2}.
 
$\bullet$ 
At times,
backstepping may not be needed but it still may be desirable to smooth
jumps in the control signal. This situation is addressed in Section
\ref{section:BS3}, where the ideal feedback is written in a form that
is affine in a function of the logic mode; the latter is then treated
as an ideal feedback and implemented dynamically through backstepping.

$\bullet$
For backstepping problems where
the SLFF pair is not ready-made, the extra dynamic variable described
in the preceding item can be exploited to achieve a backstepping result.
See Section \ref{section:BS4}.  This idea also appears
in \cite{Mayhew2011-synergistic}.
{\em Notation and terminology:}
$\Re$ ($\Re_{\geq 0}$) denotes the (nonnegative) real numbers, and
$\Re^{n}$ denotes $n$-dimensional Euclidean space. Given $x \in
\Re^{n}$, $|x|$ denotes its Euclidean norm.  The unit $n$-sphere is
$\sphere^{n} = \{x \in \Re^{n+1}: |x| = 1\}$.  
A function is called {\em smooth} if a sufficient number
of its derivatives exist and are continuous so that the derivations
make sense.  A nonnegative-valued function is said to be positive
definite with respect to a set if the function is zero if and only if
its argument belongs to the set.  For a closed set $\X \subset Q
\times \Re^{n}$, where $Q \subset \Re$ is a finite set, and a smooth
function $V:\X \rightarrow \Re$, we use $\nabla V(q,z)$ to denote
gradient of $V$ relative to $z \in \Re^{n}$, with $q \in Q$ considered
to be constant. Given a smooth function $\kappa:\X \rightarrow
\Re^{m}$, we use $\D \kappa(q,z)$ to denote the Jacobian matrix of
$\kappa$ relative to $z$, i.e., $\D \kappa(q,z)$ is an $\Re^{m \times
  n}$ matrix with $ij$-th entry given as $\partial
\kappa_{i}(q,z)/\partial z_{j}$.  As in \cite{Goebel2009}, a hybrid
system with state $x \in \Re^{n}$ is described by flow and jump sets
$C,D \subset \Re^{n}$ and set-valued flow and jump maps $F,G:\Re^{n}
\rightrightarrows \Re^{n}$. It satisfies the
basic conditions \cite{Goebel2009} if $C$ and $D$ are closed, $F$ and $G$ are outer
semicontinuous and locally bounded, $F(x)$ is nonempty and
convex for all  $x \in C$, and $G(x)$ is nonempty for all  $x \in D$.

\section{SLFF pairs}
\label{section:SLFF1}

We extend the definition of a synergistic Lyapunov function and
feedback pair defined in \cite{Mayhew2011-synergistic}.  Consider the
system
\begin{equation}
\left.
\begin{aligned}
\dot{q} &= 0 \\
\dot{z} &= f(q,z,\omega)
\end{aligned}
\right\}    \qquad (q,z) \in \X,
   \label{eq:1aa}
\end{equation}
where 
$f: \X \times \Re^{m} \to \Re^{n}$ is continuous,
$\omega \in \Re^{m}$ is the control variable, the set $\X \subset Q
\times \Re^{n}$ is closed and the set $Q \subset \Re$ finite.  Let $\A
\subset \Y \subset \X$ be such that $\A$ is compact and $\Y$ is
closed.  We define the set
\begin{equation}\label{eqn:B}
\B := \{(q,z) \in \X: \exists s \in Q, \ (s,z) \in \A\}.
\end{equation}
A 
$\mathcal{C}^{1}$ 
function $V:\X \rightarrow
\Re_{\geq 0}$ and a continuous function $\kappa:\X \rightarrow
\Re^{m}$ form a {\em synergistic Lyapunov function and feedback (SLFF)
  pair candidate} relative to $(\A,\Y)$ if
\begin{itemize}
\item
$\left\{(q,z) \in \X : V(q,z) \leq c\right\}$
is compact for each $c \geq 0$;
\item 
$V(q,z)=0$ if and only if $(q,z) \in \A$,
\item
For all $(q,z) \in \Y$,
\begin{equation}
\langle \nabla V(q,z), f(q,z,\kappa(q,z))  \rangle \leq 0  .
\label{eq:derivative}
\end{equation}
\end{itemize}
Given an SLFF pair candidate $(V,\kappa)$,
define
\begin{equation}
\label{eqn:E}
   \mathcal{E}  : =  \left\{ (q,z) \in \Y  :
   \langle \nabla V(q,z), f(q,z,\kappa(q,z)) \rangle = 0 \right\} 
\end{equation}
and let $\Psi \subset \mathcal{E}$ be the largest weakly invariant set
\cite{SanfeliceGoebelTeel05}
for
\begin{equation}
 \left.
\begin{aligned}
     \dot{q} &= 0 \\
     \dot{z} &= f(q,z,\kappa(q,z))
\end{aligned}
 \right\}   \quad (q,z) \in \mathcal{E}.
 \end{equation}
 For each $(q,z) \in \X$, define
 \begin{equation}
 \mu_{V}(q,z):= 
 \displaystyle V(q,z) - \min_{s \in Q} V(s,z)  .
 \end{equation}
  The pair $(V,\kappa)$ is called a
 {\em synergistic Lyapunov function and feedback pair relative
 to $(\A,\Y)$} 
 if
\begin{align}
 \label{eq:positive_gap}
 \mu_{V}(q,z) &> 0 & \forall (q,z) &\in 
\left( \Psi \cup \overline{ \X \setminus \Y}\right) \setminus \A,
\intertext{
 in which case $\mu_{V}(q,z)$ is called the {\em synergy gap} at $(q,z)$.
 Given a continuous function $\delta:\X \rightarrow \Re_{\geq 0}$, when}
\mu_{V}(q,z) &> \delta(q,z) & \forall (q,z) &\in 
\left( \Psi \cup \overline{\X \setminus \Y}\right) \setminus \A
 \label{eq:7}
\intertext{we say that \emph{the synergy gap exceeds $\delta$}.  When
  $\delta$ satisfies}
\label{eqn:totally}
\mu_{V}(q,z) &> \delta(q,z) & \forall (q,z) &\in
\left( \Psi \cup \overline{\X \setminus \Y} \cup \B \right) \setminus \A,
\end{align}
we say that the synergy gap \emph{totally} exceeds $\delta$.  Where
the synergy gap is positive, we can change $q$ to reduce $V$, which is desirable
at points in $\Psi \setminus \A$, where the
value of $V$ could get stuck during flows, at points in
$\left(\overline{\X \setminus \Y}\right) \setminus \A$, where the
$q$-th feedback function is not effective, and possibly at points in
$\B \setminus \A$ to ensure that the set $\B$ is stabilized.
\begin{proposition}
\label{prop:1}
{\em The synergy gap is a continuous function.  If $(V,\kappa)$ is an
  SLFF pair, then there exists a continuous function $\delta:\X
  \rightarrow \Re_{\geq 0}$ that is positive on $\X \setminus \A$ such
  that the synergy gap (totally) exceeds $\delta$.  If the synergy gap
  for $(V,\kappa)$ (totally) exceeds the function $\delta$ then, for
  each smooth $\Kinf$-function $\rho$ having a positive, nondecreasing
  derivative denoted $\rho'$, the pair $(\rho(V),\kappa)$ is an SLFF
  with synergy gap (totally) exceeding the function $\tilde{\delta}$
  defined as $\tilde{\delta}(q,z):= \rho'(c V(q,z)) (1-c)
  \delta(q,z)$, where $c$ can be taken arbitrarily in the interval
  $(0,1)$.  }
\end{proposition}

We show that the existence of an SLFF pair relative to the compact set
$\A$ is equivalent to the existence of a feedback
\begin{equation}
\begin{array}{ll}
\omega = \alpha(q,z)  &  \qquad  (q,z) \in C \subset \X \\
q^{+} \in G_{c}(q,z) \subset Q & \qquad (q,z) \in D \subset \X
\end{array}
\label{eq:controller11blue}
\end{equation}
satisfying the basic conditions and the conditions $\A \subset C$,
$C \cup D = \X$, and rendering
the compact set $\A$ globally asymptotically stable for the system
(\ref{eq:1aa}), (\ref{eq:controller11blue}), that is, for
\begin{equation}
\underbrace{
\begin{aligned}
\dot{q} &= 0 \\ \dot{z} &= f(q,z,\alpha(q,z))
\end{aligned}
}_{\displaystyle (q,z) \in C} \qquad
\underbrace{
\begin{aligned}
q^{+} &\in G_{c}(q,z) \\ z^{+} &= z
\end{aligned}
}_{\displaystyle (q,z) \in D .}
\label{eq:feedback11blue}
\end{equation}
We start by showing that this asymptotic stabilizability
property implies the existence of an SLFF.  The opposite implication
is established in Theorem \ref{theorem:1} given in Section \ref{section:HC}.

\begin{theorem}
\label{theorem:supply}
{\em Suppose the data of (\ref{eq:feedback11blue}) satisfies the basic
  conditions, the compact set $\A$ is globally asymptotically stable
  for (\ref{eq:feedback11blue}), $\A \subset C$, and $C \cup D =
  \X$. Then there exists a smooth function $V:\X \rightarrow \Re_{\geq
    0}$ such that $(V,\alpha)$ is an SLFF pair relative to $(\A,\Y)$
  with $\Y=C$ and there exists $\varepsilon>0$ such that the synergy
  gap (totally) exceeds $\delta$ where $\delta(q,z) := \varepsilon
  V(q,z)$.  If, in addition, $D \cap \A = \emptyset$ (and $\B
  \setminus \A$ is closed) then there exist $\varepsilon_{1}>0,
  \varepsilon_{2}>0$ such that the synergy gap (totally) exceeds
  $\delta$ with $\delta(q,z) :=\varepsilon_{1} V(q,z) +
  \varepsilon_{2}$.  }
\end{theorem}

\vspace{-.1in}

\section{Hybrid feedback from an SLFF pair}
\label{section:HC}

Let
$(V,\kappa)$ denote the SLFF pair and let $\delta: \X \rightarrow
\Re_{\geq 0}$ be continuous. We specify a hybrid
controller to globally asymptotically stabilize
$\A$ (and $\B$) for (\ref{eq:1aa}) as
\begin{equation}
\begin{aligned}
  C &:= \left\{ (q,z) \in \X : \mu_{V}(q,z) \leq \delta(q,z) \right\}\\
  \omega &:= \kappa(q,z) \\ 
  D &:= \left\{ (q,z) \in \X : \mu_{V}(q,z) \geq \delta(q,z)
  \right\}\\ 
  G_{c}(z) &:= \left\{ g_{c} \in Q : \mu_{V}(g_{c},z) = 0 \right\}
\end{aligned}
  \label{eq:controller_b}
\end{equation}
resulting in the closed-loop hybrid system
\begin{equation}
\underbrace{
\begin{aligned}
\dot{q} &= 0 \\
\dot{z} &= f(q,z,\kappa(q,z))
\end{aligned}}_{\displaystyle (q,z) \in C}
\qquad
\underbrace{
\begin{aligned}
q^{+} &\in G_{c}(z) \\
z^{+} &= z \\
\end{aligned}}_{\displaystyle (q,z) \in D.}
 \label{eq:closed_loop}
 \end{equation}
 Since $\delta$ and $\mu_{V}$ are continuous, $C$ and $D$ are closed. 
Since
$\mu_{V}$ is continuous, 
$G_{c}$ is outer semicontinuous.   Also,
 $C \cup D = \X$ and $G_{c}(z)$ is non-empty for each $z$ such that
 $(q,z) \in \X$ for some $q \in Q$, in particular, for $(q,z) \in D$.
\begin{theorem}
\label{theorem:1}
{\em Let $\Y \subset \X$, let $\A \subset \Y$ be compact, and let
  $\delta : \X \to \reals_{\geq 0}$ be continuous and positive on $\X
  \setminus \A$.  If $(V,\kappa)$ is an SLFF pair for (\ref{eq:1aa})
  relative to $(\A,\Y)$ with synergy gap (totally) exceeding $\delta$,
  then $\A \subset C$ and $\A$ ($\B$) is globally asymptotically
  stable for the closed-loop system (\ref{eq:closed_loop}).}
\end{theorem}

%

%

\section{Refinement of SLFF pair properties}
\label{section:SLFF2}

\subsection{Weak SLFF pairs for affine control systems}

We introduce a {\em weak} synergistic Lyapunov function and
feedback pair (weak SLFF) for 
(\ref{eq:1aa}) when 
$f(q,z,\omega) = \phi(q,z) +
\psi(q,z) \omega$ where $\phi$ and $\psi$ are smooth.  Given an SLFF
pair candidate $(V,\kappa)$, with $V$ and $\kappa$ smooth, define
\begin{equation}\label{eqn:W}
 \W := \left\{ (q,z) \in \X: \nabla V(q,z)^\top \psi(q,z) =
 0 \right\}.
\end{equation}
Recall the definition of $\E$ in Section \ref{section:SLFF1} and let
$\Omega \subset \E \cap \W$ denote the largest weakly invariant set
for the system
\begin{equation}
 \left.
\begin{aligned}
\dot{q} &= 0 \\
\dot{z} &= \phi(q,z) + \psi(q,z) \kappa(q,z)
\end{aligned}
 \right\}   \quad (q,z) \in \E \cap \W .
 \end{equation}
The pair $(V,\kappa)$ is called a {\em weak synergistic Lyapunov function and feedback pair relative to $(\A,\Y)$} if
\begin{align}
 \mu_{V}(q,z) &> 0  & \forall (q,z) 
&\in \left( \Omega \cup \overline{\X \setminus \Y} \right) \setminus \A.\\
\intertext{Given a continuous function $\delta:\X \rightarrow \Re_{\geq 0}$,
 when}
  \label{eq:42}
  \mu_{V}(q,z) &> \delta(q,z) & \forall (q,z) 
&\in \left( \Omega \cup \overline{\X \setminus \Y} \right) \setminus \A,\\
\intertext{we say that \emph{the synergy gap weakly exceeds $\delta$}. If $\delta$ satisfies}
  \label{eq:weaktotally}
  \mu_{V}(q,z) &> \delta(q,z) & \forall (q,z) 
&\in \left( \Omega \cup \overline{\X \setminus \Y} \cup \B \right) \setminus \A,
\end{align}
we say that the synergy gap weakly \emph{totally} exceeds $\delta$.
The next lemma follows immediately from the fact that $\Omega \subset
\Psi$ and then comparing (\ref{eq:42}) to (\ref{eq:7}).

\begin{lemma}
\label{lemma:1}
{\em
If $(V,\kappa)$ is a smooth SLFF pair with synergy gap (totally) exceeding
$\delta$ then it is also
a weak SLFF pair with synergy gap weakly (totally) exceeding $\delta$.
}
\end{lemma}

\begin{example}[3-D Pendulum]
\label{exm:3dpendweak}
The reduced dynamics of the 3-D pendulum are given in
\cite{Chaturvedi2010} as
\begin{subequations}
\label{eqn:3dpend}
\begin{align}
\label{eqn:3dpend-kin}
\dot{z} &= \cross{z}\omega \\
\label{eqn:3dpend-dyn}
J\dot{\omega} &= \cross{J\omega}\omega + mg\cross{\nu}z + \tau,
\end{align}
\end{subequations}
where $z \in \sphere^{2}$ is the direction of gravity in the
body-fixed frame, $\omega \in \reals^{3}$ is the angular velocity
expressed in the body-fixed frame, $m$ is the mass, $g$ is the
gravitational constant, $\nu$ is the vector from the pivot location to
the center of mass, $\tau \in \reals^{3}$ is a vector of input
torques, and for any $x,y \in \reals^{3}$, $\cross{x}$ is the $3
\times 3$ skew-symmetric matrix that satisfies $\cross{x}y = x \times
y$, where $\times$ denotes the vector cross product.
We now
stabilize the 
``inverted'' 
point
$(z,\omega) = (-\nu/|\nu|,0)$.  

Let $Q$ be a finite set, $\X_{0} = Q \times \sphere^2$, $S \subset Q$,
and $\A_{0} = S \times \{-\nu/|\nu|\}$.  Let $V_{0} : \X_{0} \to
\reals$ be positive definite on $\X_{0}$ relative to $\A_{0}$ and
define $\kappa_{0}(q,z) = 0$.  Clearly, we have that $\inner{\nabla
  V_{0}(q,z)}{\cross{z}\kappa_{0}(q,z)} = 0$ for all $(q,z) \in
\X_{0}$ so that $\Y_{0} = \X_{0} = \E_{0}$ and $\Omega_{0} = \W_{0} =
\{(q,z) \in \X_{0} : \nabla V_{0}(q,z) ^\top \cross{z} = 0 \}$. The
pair $(V_{0},\kappa_{0})$ is then a weak SLFF pair for
\eqref{eqn:3dpend-kin} relative to $(\A_{0},\X_{0})$ if
\begin{equation}\label{eqn:3dpend-syn}
\inf_{(q,z) \in \Omega_{0} \setminus \A_{0}} \mu_{V_{0}}(q,z) > 0.
\end{equation}
To satisfy \eqref{eqn:3dpend-syn}, we may use the synergistic
potential functions of \cite{Mayhew2010-3dpend, Mayhew2010-spherical}.
We henceforth assume that the synergy gap weakly totally exceeds a
constant $\delta(q,z) = c > 0$.
\null \hfill \null $\square$
\end{example}

%

\subsection{Pure and Ready-made SLFF pairs}

When
$\Y = \X$, a (weak) SLFF pair is called a (weak)
{\em pure} SLFF pair.  
A weak SLFF pair $(V,\kappa)$ with synergy gap weakly (totally)
exceeding $\delta:\X \rightarrow \Re_{\geq 0}$
is said to be {\em type I ready-made} with
respect to the continuous, positive definite
 function $\sigma:\Re^{m} \rightarrow \Re_{\geq 0}$
if there exists a continuous function
$\varrho:\X \rightarrow \Re_{\geq 0}$ such that,
 $\forall (q,z) \in \X$,
$s \in Q$,
and $\omega=\kappa(q,z)$,
\begin{equation}
\sigma \bigl( \omega - \kappa(s,z) \bigr) -
 \sigma \bigl( \omega - \kappa(q,z) \bigr)
 \leq \varrho(q,z) 
 \label{eq:for_ready1}
\end{equation}
and, for all $(q,z) \in \left( \Omega \setminus \A \right) \cup
\overline{\X \setminus \Y}$, 
\begin{equation}
\mu_{V}(q,z) > \delta(q,z)+ \varrho(q,z) .
\label{eq:for_ready2}
\end{equation}
Since $\mu_{V}(q,z)=0$ for $(q,z) \in \A$, the type I ready-made property
implies that 
\begin{equation}
\overline{\X \setminus \Y} \cap \A = \emptyset .
\label{eq:necessarily}
\end{equation}
If $\kappa$ does not depend on $q$ then,
in (\ref{eq:for_ready1}), we
can take $\varrho(q,z)=0$ for all $(q,z) \in \X$.
With this choice for $\varrho$, if (\ref{eq:necessarily}) holds then (\ref{eq:for_ready2})
follows from (\ref{eq:42}).
According to the last statement of Proposition \ref{prop:1},
if $\delta$ is positive valued, $V$
is radially unbounded, and the condition (\ref{eq:necessarily}) holds, then
the type I ready-made property
is achievable for any $\sigma$ by modifying the function $V$ as
$\rho(V)$ with $\rho$ chosen appropriately.

A weak SLFF pair $(V,\kappa)$ with synergy gap weakly (totally)
exceeding $\delta:\X \rightarrow \Re_{\geq 0}$ is said to be {\em type
  II ready-made} with respect to the continuous, positive definite
function $\sigma:\Re^{m} \rightarrow \Re_{\geq 0}$ if there exists a
continuous function $\varrho:\X \rightarrow \Re_{\geq 0}$ such that,
for all $(q,z) \in \X$, $s \in Q$, and $\omega \in \Re^{m}$,
(\ref{eq:for_ready1}) holds and, moreover, (\ref{eq:for_ready2}) holds
for all $(q,z) \in \left( \Omega \setminus \A \right) \cup
\overline{\X \setminus \Y}$.  In particular, the difference between
type I and type II ready-made is in the requirement on $\omega$ for
which (\ref{eq:for_ready1}) holds: $\omega=\kappa(q,z)$ for type I and
$\omega \in \Re^{m}$ for type II.  Clearly, if the SLFF pair is type
II ready-made then it is type I ready-made.  Like for the type I case,
if $\kappa$ is independent of $q$ then, in (\ref{eq:for_ready1}), we
can take $\varrho(q,z)=0$ for all $(q,z) \in \X$.

\begin{example}[3-D pendulum]
The weak SLFF pair $(V_{0},\kappa_{0})$ for \eqref{eqn:3dpend-kin}
with synergy gap weakly totally exceeding $c > 0$ is type I/II ready
made with respect to any positive definite function $\sigma$ and
appropriate function $\varrho$, since $\kappa_{0}(q,z) \equiv 0$ does
not depend on $q$.
\null \hfill \null $\square$
\end{example}

\section{Ready-made backstepping}
\label{section:BS0}

The ensuing backstepping results are useful mainly for the case where
the SLFF pair for the reduced-order system has a synergy gap (totally)
exceeding a positive-valued continuous function $\delta$, i.e.,
$\delta: \X \rightarrow \Re_{>0}$.  Indeed, the nature of our
backstepping results is that the extended system admits an SLFF pair
with synergy gap (totally) exceeding the same function $\delta$.  If
$\delta$ is not positive valued then, since it does not depend on the
extended state, there is no hope of it being positive valued away from
the attractor in the extended state space.  In this case, the hybrid
control construction based on an SLFF pair given in Theorem
\ref{theorem:1} would not be applicable.

\subsection{From a weak, pure, ready-made SLFF pair}
\label{section:BS1}

We consider the control system
\begin{equation}
\left.
\begin{array}{l}
\dot{q} = 0 \\
\dot{\zeta} = \phi_{1}(q,\zeta)+\psi_{1}(q,\zeta)u
\end{array}
\right\} \qquad (q,\zeta) \in \X_{1}
\label{eq:15blue}
\end{equation}
with $u \in \Re^{m}$,
where 
$\zeta=(z^{\top},\omega^{\top})^{\top}$, $\X_{1} = \X_{0} \times \Re^{m}$ and
\begin{equation}
 \phi_{1}(q,\zeta) = 
\begin{bmatrix}
\phi_{0}(q,z) + \psi_{0}(q,z) \omega \\
0
\end{bmatrix},
\
\psi_{1}(q,\zeta) = 
\begin{bmatrix} 0 \\ 1  
\end{bmatrix}.
\label{eq:16blue}
\end{equation}
We construct a (non-weak) SLFF pair with synergy gap exceeding a
positive-valued function $\delta$ by supposing we have a weak, pure,
ready-made SLFF pair with synergy gap weakly (totally) exceeding
$\delta$ for the reduced system
\begin{equation}
\left.
\begin{array}{l}
  \dot{q} = 0 \\
  \dot{z} = \phi_{0}(q,z) + \psi_{0}(q,z) \omega
\end{array}
\right\} \qquad (q,z) \in \X_0
\label{eq:1aaaa_blue}
\end{equation}
with controls $\omega \in \Re^{m}$.

Let $\A_{0} \subset \X_0$ be compact.  For the system
(\ref{eq:1aaaa_blue}), let $(V_{0},\kappa_{0})$ be a weak SLFF pair
relative to $(\A_{0},\X_{0})$, with synergy gap weakly (totally)
exceeding the continuous function $\delta:\X_{0} \rightarrow \Re_{
  \geq 0}$.  Let $\Gamma \in \Re^{m \times m}$ be a symmetric,
positive definite matrix and suppose that the SLFF pair is type I
ready-made relative to $\sigma(v):=v^{\top} \Gamma v$.  Define
 \begin{equation}
    \mathcal{A}_{1} := \left\{(q,\zeta) \in \X_1:  (q,z) \in \mathcal{A}_{0},  \
     \omega = \kappa_{0}(q,z) \right\} .
    \label{eq:A1blue}
 \end{equation}
 
 For each $(q,\zeta) \in \X_{1}$, define
 \begin{equation}
  V_{1}(q,\zeta) = V_{0}(q,z) + \sigma(\omega - \kappa_{0}(q,z)) .
  \label{eq:V1blue}
  \end{equation}
  
Let $\theta:\Re_{\geq 0} \rightarrow \Re_{\geq 0}$ 
be a continuous, positive definite function, and 
let the smooth function $\Theta:\Re^{m} \rightarrow \Re^{m}$ satisfy
\begin{equation}
 v^{\top} \Gamma \Theta(v) + \Theta(v)^{\top} \Gamma v \leq - \theta(|v|)  \quad
 \forall v \in \Re^{m} .
\label{eq:damping}
\end{equation}
Define
 \begin{equation}
\begin{aligned}
  \kappa_{1}(q,\zeta) &:=    \Theta (\omega - \kappa_{0}(q,z) ) \\
  &\quad -\tfrac{1}{2}\Gamma^{-1}\psi_{0}(q,z)^{\top}\nabla V_{0}(q,z)  \\
  &\quad + \D \kappa_{0}(q,z)  (\phi_{0}(q,z) + \psi_{0}(q,z) \omega)  .
\end{aligned}
  \label{eq:kappa1blue}
\end{equation}
 
\begin{theorem}
\label{theorem:2blue}
{\em Let the compact set $\mathcal{A}_{0}$ and the smooth functions
  $(V_{0},\kappa_{0})$ be given.  Let the compact set
  $\mathcal{A}_{1}$ be defined as in (\ref{eq:A1blue}) and let the
  pair $(V_{1},\kappa_{1})$ be defined by
  (\ref{eq:V1blue})-(\ref{eq:kappa1blue}).  Suppose, for the system
  (\ref{eq:1aaaa_blue}), that the pair $(V_{0},\kappa_{0})$ is a weak
  SLFF relative to the pair $(\A_{0},\X_{0})$ with synergy gap weakly
  (totally) exceeding the continuous function $\delta:\X_{0}
  \rightarrow \Re_{\geq 0}$ and the SLFF pair is type I ready-made
  relative to the function $\sigma$ defined as $\sigma(v) = v^{\top}
  \Gamma v$ where $\Gamma$ is a symmetric positive definite matrix.
  Under these conditions, for the system
  (\ref{eq:15blue})-(\ref{eq:16blue}), the pair $(V_{1},\kappa_{1})$
  is an (non-weak) SLFF pair relative to the pair $(\A_{1},\X_{1})$
  with synergy gap (totally) exceeding $\delta$.  }
\end{theorem}

\begin{example}[3-D pendulum]
\label{exm:3dpend-bs}
Consider the input transformation $\tau = -\cross{J\omega}\omega
-mg\cross{\nu}z + Ju$, which renders the angular velocity dynamics
\eqref{eqn:3dpend-dyn} as $\dot{\omega} = u$.  We now apply
Theorem~\ref{theorem:2blue}.  Let $\sigma(\omega) =
\frac{1}{2}\omega^{\top} J \omega$ (i.e. $\Gamma = J/2$) and let
$\Theta(\omega) = J^{-1}\left(\cross{J\omega}\omega -
\Xi(\omega)\right)$, where $\Xi : \reals^{3} \to \reals^{3}$ satisfies
$\omega^{\top} \Xi(\omega) \geq \theta(|\omega|)$ and $\theta: \reals
\to \reals$ is a continuous, positive definite function.  Applying
\eqref{eq:kappa1blue}, we arrive at
\[
u(q,z) = J^{-1}\left(\cross{J\omega}\omega - \Xi(\omega)\right) -
J^{-1}\cross{z}^{\top}\nabla V(q,z),
\]
which yields
\begin{equation}\label{eqn:3dpend-torquefb}
\tau = \kappa_{1}(q,z) = -mg\cross{\nu}z - \Xi(\omega) - \cross{z}^\top \nabla
V_{0}(q,z).
\end{equation}
This recovers the feedback of \cite{Mayhew2010-3dpend}.  As a result
of Theorem~\ref{theorem:2blue}, it follows that $(V_{1},\kappa_{1})$,
with $V_{1}(q,z,\omega) = V_{0}(q,z) + \frac{1}{2}\omega^\top J
\omega$, is an SLFF pair for \eqref{eqn:3dpend} relative to
$(\A_{1},\X_{1})$, where $\X_{1} = Q \times \sphere^{2} \times
\reals^{3}$ and $\A_{1} = \{(q,z,\omega) \in \X_1: q \in S, \ z =
-\nu/|\nu|, \ \omega = 0\}$, with gap totally exceeding $\delta(q,z) =
c$.
\null \hfill \null $\square$
\end{example}

\subsection{From a weak, ready-made SLFF pair}
\label{section:BS2}

We again consider the control system (\ref{eq:15blue})-(\ref{eq:16blue}).
Let $\A_{0} \subset \X_{0}$ be compact and let $\Y_{0} \subset \X_{0}$ be closed.
The results in this section apply to the case where $\Y_{0}$ is not
necessarily equal to $\X_{0}$.  For the system (\ref{eq:1aaaa_blue}),
let $(V_{0},\kappa_{0})$ be a weak, SLFF pair relative to
$(\A_{0},\Y_{0})$ with synergy gap weakly (totally) exceeding the
continuous function $\delta:\X_{0} \rightarrow \Re_{ \geq 0}$.  In
addition, suppose the SLFF pair is type I ready-made relative to
$\sigma(v):=L|v|$ where $L>0$.  Let $\Gamma \in \Re^{m \times m}$ be a
positive definite, symmetric matrix.  Define $\sigma_{2}(v) :=
v^{\top} \Gamma v$ for all $v \in \Re^{m}$.  Let $\rho \in \Kinf$ be
smooth, such that $\rho'(s)>0$ for all $s \geq 0$, and such that $\rho
\circ \sigma_{2}$ is globally Lipschitz with constant less than or
equal to $L$.  For example, pick $\rho(s) = c \tilde{\rho}(s)$ where
$c>0$ is sufficiently small and $\tilde{\rho}(s)=s$ for $s \in [0,1]$,
$\tilde{\rho}(s) = k \sqrt{s}$ for $s \geq 2$ where $k \geq 1$, and
such that $\tilde{\rho}'(s)>0$ for $s \in [1,2]$ to smoothly connect
the value $1$ at $s=1$ to the value $k \sqrt{2}$ at the value $s=2$.
This construction makes the SLFF pair $(V_{0},\kappa_{0})$ type II
ready-made for the function $v \mapsto \rho (\sigma_{2}(v))$.  Define
$\Y_{1} := \Y_{0} \times \Re^{m}$,
\begin{equation}
 V_{1}(q,\zeta) := V_{0}(q,z) + \rho(\sigma_{2}(\omega-\kappa_{0}(q,z))) .
 \label{eq:V1red}
\end{equation}
and
 \begin{equation}
\begin{aligned}
  \kappa_{1}(q,\zeta) &:=    \displaystyle
  \frac{1}{\rho'(\sigma_{2}(\omega-\kappa_{0}(q,z)))}
  \left[ \vphantom{\frac{1}{2}}  \Theta (\omega - \kappa_{0}(q,z) ) \right. \\
  &\quad -  \tfrac{1}{2}\Gamma^{-1} \psi_{0}(q,z)^{\top} \nabla V_0(q,z)  \\
  &\quad \displaystyle \left. + \mathcal{D} \kappa_0(q,z)  (\phi_{0}(q,z) + \psi_{0}(q,z) \omega) 
  \vphantom{\frac{1}{2}} \right] .
\end{aligned}
  \label{eq:kappa1red}
\end{equation}

\begin{theorem}
\label{theorem:2red}
{\em
Let the compact set $\A_{0} \subset \X_{0}$ and the
closed set $\Y_{0} \subset \X_{0}$ be given. Let $\rho$ and $\sigma_{2}$
be such that $\sigma_{2}(v) = v^{\top} \Gamma v$ for all $v \in \Re^{m}$
where $\Gamma \in \Re^{m \times m}$ is a symmetric, positive definite matrix,
$\rho \in \Kinf$ is smooth, $\rho'(s)>0$ for all $s \geq 0$, and
$v \mapsto \rho(\sigma_{2}(v))$ is globally
Lipschitz with constant less than or equal to $L>0$.
Let the compact set $\mathcal{A}_{1}$ be defined as in (\ref{eq:A1blue}) and let the pair $(V_{1},\kappa_{1})$ be
defined by (\ref{eq:V1red})-(\ref{eq:kappa1red}).
Suppose, for the system (\ref{eq:1aaaa_blue}), that
the pair $(V_{0},\kappa_{0})$ is a weak SLFF
pair relative to $(\A_{0},\Y_{0})$ 
with synergy gap weakly (totally) exceeding the continuous
function $\delta:\X_{0} \rightarrow \Re_{\geq 0}$ and the SLFF pair
is type I ready-made relative to the function $\sigma$ given by $\sigma(v) := L |v|$
for all $v \in \Re^{m}$.
Under these conditions,
for the system (\ref{eq:15blue})-(\ref{eq:16blue}),
 the pair $(V_{1},\kappa_{1})$
is a (non-weak) SLFF pair relative to $(\A_{1},\Y_{1})$
with synergy gap (totally) exceeding $\delta$.
}
\end{theorem}

\section{Smoothing without backstepping}
\label{section:BS3}

Now we consider the situation where the control does not enter through
an integrator but we want to remove jumps from the feedback.  The
ideas described here are also used in Section \ref{section:BS4} for a
backstepping algorithm that does not require the SLFF pair to be
ready-made.  Henceforth,
we work with SLFF pairs having a synergy gap bounded away from
a function $\delta$.
The synergy
gap is said to be (totally) bounded away from a continuous function
$\delta:\X \rightarrow \Re_{\geq 0}$ if there exists a positive
real number $\varepsilon$ such that the energy gap (totally) exceeds the
function $(q,z) \mapsto \tilde{\delta}(q,z):=\delta(q,z) + \varepsilon$.
We note that if the synergy gap is (totally) bounded away from a continuous
function $\delta:\X \rightarrow \Re_{\geq 0}$ then, because $\mu_{V}(q,z)=0$
for $(q,z) \in \A$, it follows that $\overline{\X \setminus \Y} \cap \A = \emptyset$.
We start with the control system
\begin{equation}
\left.
\begin{array}{l}
  \dot{q} = 0 \\
  \dot{z} = \phi_{0}(q,z) + \psi_{0}(q,z) \omega
\end{array}
\right\} \qquad (q,z) \in \X_0
\label{eq:1aaaa_green}
\end{equation}
with controls $\omega \in \Re^{m}$
for which we suppose we have a (non-weak) SLFF
pair relative to the pair ($\A_{0},\X_{0})$ where $\A_{0} \subset \X_{0}$
is compact with synergy gap (totally) bounded away from the function
$\delta:\X \rightarrow \Re_{\geq 0}$.   Let
$M_{0}$ be the projection of $\X_{0}$ in the $z$ direction, i.e.,
\begin{displaymath}
M_{0}:= \left\{ z \in \Re^{n} : (q,z) \in \X_{0} \ \mbox{\rm for some} \ q \in Q \right\} .
\end{displaymath}
Let $N$ be the cardinality of $Q$, let $r \leq N$, and let the smooth
functions $\beta_{0}:M_{0} \rightarrow \Re^{m}$ and
$\vartheta_{0}:M_{0} \rightarrow \Re^{m \times r}$ satisfy 
\begin{equation}\label{eqn:affinelogic}
\kappa_{0}(q,z)= \beta_{0}(z) + \vartheta_{0}(z) \varsigma_{0}(q) 
\qquad \forall (q,z) \in \X_{0},
\end{equation}
where $\varsigma_{0}: Q \to \reals^{r}$ is some function of the
variable $q$.  In turn, we see that, for the system
\begin{equation}\nonumber
\left.
\begin{array}{l}
  \dot{q} = 0 \\
  \dot{z} = \phi_{0}(q,z) + \psi_{0}(q,z)  \left( \beta_{0}(z) + \vartheta_{0}(z) p \right)
\end{array}
\right\} \quad (q,z) \in \X_0
\label{eq:1aazz_green}
\end{equation}
with $p$ as the control variable, the pair $(V_{0},\varsigma_{0})$ is
an SLFF pair with respect to $(\A_{0},\X_{0})$ with synergy gap
(totally) bounded away from $\delta$.  Let $\varepsilon>0$ be such
that the synergy gap (totally) exceeds
$\tilde{\delta}_{1}(q,z):=\delta(q,z)+ \varepsilon$.  Since the set
$Q$ is finite, we can easily find a positive definite, symmetric
matrix $\Gamma$ such that, with $\sigma_{0}(v) = v^{\top} \Gamma v$,
we have
\begin{equation}
\sigma_{0}(\varsigma_{0}(q) - \varsigma_{0}(s)) \leq \frac{\varepsilon}{2} 
\qquad \forall (q,s) \in Q \times Q .
\end{equation}
This implies that the SLFF pair $(V_{0},\varsigma_{0})$, with a
synergy gap (totally) exceeding $\tilde{\delta}_{2}(q,z):=\delta(q,z) +
\varepsilon/2$, is type I ready-made for backstepping relative to
$\sigma_{0}$.  Like in Section \ref{section:BS2}, we can also find a
function $\rho$ so that, for all $(p,q,s) \in \Re^{r} \times Q \times
Q$,
\begin{equation}
\rho (\sigma_{0}(p - \varsigma_{0}(s)))
-\rho (\sigma_{0}(p - \varsigma_{0}(q))) \leq \frac{\varepsilon}{2} 
.
\end{equation}
In particular, this implies
the SLFF pair $(V_{0},\varsigma_{0})$, with a synergy gap (totally) exceeding
$\tilde{\delta}_{2}(q,z):=\delta(q,z) + \varepsilon/2$,
 is type II ready-made for backstepping relative
to $\rho \circ \sigma_{0}$.

Now, using Lemma \ref{lemma:1}, and depending on whether the original pair $(V_{0},\kappa_{0})$
was pure or not, we can
apply either Theorem \ref{theorem:2blue} or Theorem \ref{theorem:2red}
 to construct a pair $(V_{1},\varsigma_{1})$ that is
an SLFF pair with synergy gap (totally) exceeding $\tilde{\delta}_{2}$ for
\begin{displaymath}
\begin{array}{l}
\left.
\begin{array}{l}
  \dot{q} = 0 \\
  \dot{z} = \phi_{0}(q,z) + \psi_{0}(q,z) \left( \beta_{0}(z) +\vartheta_{0}(z) p \right), \ 
  \dot{p} = \alpha
\end{array}
\right\}  \\
\qquad \qquad \qquad \qquad \qquad \qquad \qquad
(q,z,p) \in \X_0 \times \Re^{r} .
\end{array}
\end{displaymath}
In particular, from the definition of $\tilde{\delta}_{2}$, it follows that
the synergy gap is (totally) bounded away from $\delta$.

Note that if 
$(V_{0},\kappa_{0})$ was a {\em weak} SLFF pair for the system 
(\ref{eq:1aaaa_green}),
this fact would not necessarily guarantee that
$(V_{0},\varsigma_{0})$ is  a {\em weak} SLFF pair for (\ref{eq:1aazz_green}),
because of the $\vartheta_{0}$ term that multiplies $\psi_{0}$ to generate
the input vector field.
This observation motivates assuming that $(V_{0},\kappa_{0})$ is an (non-weak)
SLFF pair for the system (\ref{eq:1aaaa_green}).   In the next section,
we will want to allow $(V_{0},\kappa_{0})$ to be a weak SLFF pair
for the system (\ref{eq:1aaaa_green}) in anticipation of another backstepping
result.  We will be able to get away with this weakened assumption because
 we will come back to the integral of $\omega$, rather than the
integral of $p$, as the control variable.

\begin{example}[3-D pendulum]
\label{exm:3dpendsmooth}
Let $\bfe_{i} \in \reals^{N}$ denote the vector with $1$ in the $i$th
index and zeros elsewhere.  Assuming that (without loss of generality)
$Q = \{1,\dotsc,N\}$, $\kappa_{1}$---as defined in
\eqref{eqn:3dpend-torquefb}---can be written as
\eqref{eqn:affinelogic}.  In particular, define
\[
\begin{aligned}
\cV_{0}(z) &= [V_{0}(1,z) \ \cdots \ V_{0}(N,z)]^\top &
\vartheta(z) &= \cross{z}\D\cV_{0}(z)^{\top} \\
\beta(z) &= -mg\cross{\nu}z - \Xi(\omega) &
\varsigma(q) &= \bfe_{q},
\end{aligned}
\]
which yields the closed-loop dynamics of \eqref{eqn:3dpend} as
\[
\begin{aligned}
\dot{z} &= \cross{z}\omega &
 J\dot{\omega} &= \cross{J\omega}\omega - \Xi(\omega)
+ \vartheta(z)\bfe_{q}.
\end{aligned}
\]
By replacing $\bfe_{q}$ with a control variable $p$, we have that
$(V_{1},\varsigma)$ is a (non-weak) SLFF pair relative to
$(\A_{1},\X_{1})$ (with $V_{1}$, $\A_{1}$ and $\X_{1}$ defined in
Example~\ref{exm:3dpend-bs}) with synergy gap totally exceeding
$\delta(q,z,\omega) = c$.  Suppose also that the synergy gap totally
exceeds $c + \epsilon$ and let $\sigma(v) = \frac{\epsilon}{8} |v|^2$
so that for all $(q,s) \in Q \times Q$, $\sigma(\bfe_{q} - \bfe_{s})
\leq \epsilon/2$ and $(V_{1},\varsigma)$ is also type I ready-made
with respect to $\sigma$.

Now, define $V_{2}(q,z,\omega,p) = V_{1}(q,z,\omega) +
\sigma(p-\bfe_{q})$, $\X_{2} = Q \times \sphere^{2} \times
\reals^{3} \times \reals^{N}$, $\A_{2} = \{(q,z,\omega,p) \in \X: q
\in S, \ z = -\nu/|\nu|, \ \omega = 0 \ p = \bfe_{q}\}$, and
\[
\begin{aligned}
\gamma(q,z,\omega,p) &= \Theta(p-\bfe_{q}) - \D\cV_{0}(z) \cross{z}\omega,
\end{aligned}
\]
where $\Theta : \reals^{N} \to \reals^{N}$ satisfies
\eqref{eq:damping} with $\Gamma = I$.  It follows from
Theorem~\ref{theorem:2blue} that $(V_{2},\gamma)$ is an SLFF pair
relative to $(\A_{2},\X_{2})$ with synergy gap totally exceeding $c +
\epsilon/2$ for the system
\[
\begin{aligned}
\dot{q} &= 0 & \dot{z} &= \cross{z}\omega \\
\dot{p} &= \alpha & J\dot{\omega} &= \cross{J\omega}\omega -
\Xi(\omega) + \vartheta(z)p
\end{aligned}
\]
with $\alpha$ as the control variable.

Having input $(z,\omega) \in \sphere^{2} \times \reals^{3}$, memory
states $(q,p) \in Q \times \reals^{N}$, and output $\tau$, the hybrid
controller for the 3-D pendulum with smoothing is given as
\[
\underbrace{
\begin{aligned}
\tau &= \beta(z) + \vartheta(z)p, \quad
\dot{q} = 0 \\
\dot{p} &= \Theta(p-\bfe_{q}) - \D\cV_{0}(z)\cross{z}\omega
\end{aligned}
}_{\displaystyle (q,z,\omega,p) \in C}
\qquad
\underbrace{
\begin{aligned}
q^{+} &= G(z,\omega,p) \\
p^{+} &= p
\end{aligned}
}_{\displaystyle (q,z,\omega,p) \in D,}
\]
where
\[
\begin{aligned}
C &= \{(q,z,\omega,p) \in \X : \mu_{V_{2}}(q,z,\omega,p) \leq c +
\epsilon/2\} \\
D &= \{(q,z,\omega,p) \in \X : \mu_{V_{2}}(q,z,\omega,p) \geq c +
\epsilon/2\} \\
G(z,\omega,p) &= \{g \in Q : \mu_{V_{2}}(g,z,\omega,p) = 0\}.
\end{aligned}
\]
If $V_{0}$ satisfies \eqref{eqn:3dpend-syn}, this
controller globally asymptotically stabilizes $\B_{2}$, where
$\B_{2}$ is related to $\A_{2}$ through \eqref{eqn:B}.
\null \hfill \null $\square$
\end{example}

\section{Backstepping without being ready-made}
\label{section:BS4}

While the backstepping constructions in this section 
use extra dynamic states, their advantage is that no preliminary step is needed to
make them ready-made for backstepping. 
Suppose we have a non-weak SLFF pair $(V_{0},\kappa_{0})$ with
synergy gap (totally) bounded away from $\delta$ for
\begin{equation}
\begin{array}{l}
\left.
\begin{array}{l}
  \dot{q} = 0 \\
  \dot{z} = \phi_{0}(q,z) + \psi_{0}(q,z) \omega
\end{array}
\right\} 
\quad 
(q,z) \in \X_0.
\end{array}
\label{eq:start_here_final}
\end{equation}
\null From the results of Section \ref{section:BS3},
the pair $(V_{1},\kappa_{1})$,
of the form
\begin{equation}
\begin{array}{l}
V_{1}(q,z,p) = V_{0}(q,z) + \sigma(p - \varsigma_{0}(q)) \\
\kappa_{1}(q,z,p) = \beta_{0}(z) + \vartheta_{0}(z) p,
\end{array}
\label{eq:first_pair}
\end{equation}
is a non-weak SLFF with synergy gap (totally) bounded away from $\delta$ for the
system
\begin{equation}
\begin{array}{l}
\left.
\begin{array}{l}
  \dot{q} = 0 \\
  \dot{z} = \phi_{0}(q,z) + \psi_{0}(q,z) \omega \\
  \dot{p} = \varsigma_{1}(q,z,p)
\end{array}
\right\} 
\quad 
(q,z,p) \in \X_0 \times \Re^{r} .
\end{array}
\label{eq:intermediate_final}
\end{equation}
Moreover, the pair $(V_{1},\kappa_{1})$ is both type I and type II ready-made
with respect to any function.   Indeed, since $\kappa_{1}$ does not depend
on $q$, we can take $\varrho(q,z,p)=0$ for all $(q,z,p) \in \X_{0} \times \Re^{r}$
in (\ref{eq:for_ready1}) and then, since (\ref{eq:necessarily}) holds because
the synergy gap is (totally) bounded away from $\delta$, (\ref{eq:for_ready2}) holds.
Now we can apply Theorem \ref{theorem:2blue} or,
if the SLFF pair is not pure, Theorem \ref{theorem:2red} to generate
a non-weak SLFF pair $(V_{2},\kappa_{2})$ with synergy gap (totally)
bounded away from $\delta$ for the extended system
\begin{equation}
\begin{array}{l}
\left.
\begin{array}{l}
  \dot{q} = 0 \\
  \dot{z} = \phi_{0}(q,z) + \psi_{0}(q,z) \omega \\
  \dot{p} = \varsigma_{1}(q,z,p) \\
  \dot{\omega} = u 
\end{array}
\right\} 
\\
\qquad \qquad  \qquad 
\qquad \qquad 
(q,z,p,\omega) \in \X_0 \times \Re^{r} \times \Re^{m}.
\end{array}
\label{eq:extended_final}
\end{equation}

Finally, consider the case where $(V_{0},\kappa_{0})$
is a weak (rather than non-weak) SLFF pair for (\ref{eq:start_here_final}).
In this case it turns out that the SLFF pair
$(V_{1},\kappa_{1})$ of the form (\ref{eq:first_pair}) is a weak SLFF pair
for the system (\ref{eq:intermediate_final}).   This fact is explained below.
From here, Theorem
\ref{theorem:2blue} or \ref{theorem:2red} can be applied as above to
derive a {\em non-weak} SLFF pair $(V_{2},\kappa_{2})$ for the system
(\ref{eq:extended_final}).

Suppose $(V_{0},\kappa_{0})$ is a weak SLFF pair for (\ref{eq:start_here_final}).
Write the system (\ref{eq:intermediate_final}) in the form
\begin{equation}
\begin{array}{l}
\left.
\begin{array}{l}
  \dot{q} = 0 \\
  \dot{\zeta} = \phi_{1}(q,\zeta) + \psi_{1}(q,\zeta) \omega
\end{array}
\right\} 
\quad
(q,\zeta) \in \X_1
\end{array}
\label{eq:extended_final2}
\end{equation}
where $\zeta:=(z^{\top},p^{\top})^{\top}$, $\X_{1} := \X_{0} \times \Re^{r}$,
\begin{equation}
 \phi_{1}(q,\zeta) := 
\begin{bmatrix} \phi_{0}(q,z) \\ \varsigma_{1}(q,z,p) \end{bmatrix}, \
 \psi_{1}(q,\zeta): = 
\begin{bmatrix}
\psi_{0}(q,z) \\  0 
\end{bmatrix}.
\end{equation}
It follows from the definitions that
\begin{displaymath}
 \nabla V_{1}(q,\zeta)^{\top} \psi_{1}(q,\zeta) =
 \nabla V_{0}(q,z)^{\top}\psi_{0}(q,z).
\end{displaymath}
Also, it follows from the proof of Theorems \ref{theorem:2blue}
and \ref{theorem:2red} that
\begin{displaymath}
\begin{array}{l}
\langle \nabla V_{1}(q,z) , \phi_{1}(q,\zeta)+\psi_{1}(q,\zeta) \kappa_{1}(q,\zeta)
\rangle = 0 \\
\Longrightarrow
\left\{
\begin{array}{l}
0 = \langle \nabla V_{0}(q,z), \phi_{0}(q,z) + \psi_{0}(q,z) \kappa_{0}(q,z) \rangle \\
p = \varsigma_{0}(q).
\end{array}
\right.
\end{array}
\end{displaymath}
Therefore
 $\Omega_{1} = \left\{(q,\zeta) \in \X_{1} : (q,z) \in \Omega_{0}, \ 
p = \varsigma_{0}(q)
 \right\}.$
This relationship can be used 
to arrive at the conclusion that $(V_{1},\kappa_{1})$ is a {\em weak}
SLFF for the system (\ref{eq:intermediate_final}) with synergy gap
(totally) bounded away from $\delta$.

%

\bibliographystyle{IEEEtran}
\bibliography{cgmBiblio,mypapers,long,Biblio,RGS}


\section{Appendices}

\subsection{Proof of Proposition \ref{prop:1}}

The continuity of the synergy gap follows from the continuity of $V$.
Since $\Psi$ is closed \cite[Lemma 3.3]{SanfeliceGoebelTeel05}, it is possible to find a function $\delta$ that
is positive on $\X \setminus \cA$ so that the synergy gap exceeds
$\delta$.  A possible function $\delta: \X \rightarrow \Re_{\geq 0}$
is given as
\begin{displaymath}
 \delta(q,z) : = \inf_{(s,\zeta) \in \Psi \cup \overline{\X \setminus \Y}}
 \left( \vphantom{\frac{1}{2}} |(q,z) - (s,\zeta)| + 0.5 \mu_{V}(s,\zeta) \right)
 \end{displaymath}
 which is continuous, 
 for all $(q,z) \in \left( \Psi \cup \overline{\X \setminus \Y} \right) \setminus \A$
 satisfies
 $\delta(q,z) \leq 0.5 \mu_{V}(q,z) < \mu_{V}(q,z)$,
 and, using (\ref{eq:positive_gap}), the fact that $\Psi$ is closed, and the continuity
 of $\mu_{V}$,
 satisfies $\delta(q,z)>0$ if $(q,z) \notin \A$.  Let $W(q,z)=\rho(V(q,z))$
 where $\rho$ is a smooth $\Kinf$ function and $\rho'$ is positive and nondecreasing.   It follows that $(W,\kappa)$ is an SLFF and the sets
 $\mathcal{E}$ and $\Psi$ are the same as those for $(V,\kappa)$.   Moreover,
 \begin{equation}
 \mu_{W}(q,z) = \rho(V(q,z)) - \min_{s \in Q} \rho(V(s,z)) .
 \end{equation}
 Inspired by the calculations in 
 \cite{SontagTeel95,NesicTeel01bb},
 we lower bound $\mu_{W}(q,z)$ 
 on the set $\left( \Psi \cup \overline{\X \setminus \Y} \right) \setminus \A$ by considering two cases:  
 $\min_{s \in Q} V(s,z) \leq c V(q,z)$ and $\min_{s \in Q} V(s,z)
 \geq c V(q,z)$ where $c \in (0,1)$.   In the first case, using the mean-value
 theorem and the monotonicity of $\rho'$,
 \begin{displaymath}
   \begin{array}{ccl}
  \mu_{W}(q,z) & \geq & \rho(V(q,z)) - \rho(c V(q,z)) \\
    & \geq & \rho'(c V(q,z)) (1-c) V(q,z) \\
    & \geq & \displaystyle \rho'(c V(q,z)) (1-c) \left( V(q,z) - \min_{s \in Q} V(s,z) \right) \\
    & > & \rho'(c V(q,z)) (1-c) \delta(q,z) .
    \end{array}
\end{displaymath}
In the second case, using the mean-value theorem and the monotonicity of $\rho'$,
\begin{displaymath}
   \begin{array}{ccl}
  \mu_{W}(q,z)
    & \geq & \displaystyle \rho'(c V(q,z)) \left( V(q,z) - \min_{s \in Q} V(s,z) \right) \\
    & > & \rho'(c V(q,z)) \delta(q,z) \\
    & > & \rho'(c V(q,z)) (1-c) \delta(q,z) .
    \end{array}
\end{displaymath}
These bounds establish the final statement of the proposition.
\null \hfill \null $\blacksquare$

\subsection{Proof of Theorem \ref{theorem:supply}}

According to \cite[Theorem~4.2]{Cai2008a}, there exists a smooth function
$V:\Re^{n+1} \rightarrow \Re_{\geq 0}$ that is radially unbounded and
positive definite with respect to the compact set $\A$ and such that,
for all $(q,z) \in C=\Y$,
\begin{equation}
\langle \nabla V(q,z) , f(q,z, \alpha(q,z)) \rangle
\leq - V(q,z) \leq 0 ,
\label{eq:in_converse1blue}
\end{equation}
and, for all $(q,z) \in D$ and $s \in G_{c}(q,z)$,
\begin{equation}
   V(s,z) \leq e^{-1} V(q,z)  .
   \label{eq:in_converse2blue}
\end{equation}
The properties of $V$ together with (\ref{eq:in_converse1blue}) 
make $(V,\alpha)$ an SLFF pair candidate relative to $(\A,\Y)$.  In addition, (\ref{eq:in_converse1blue})
guarantees
that the set $\mathcal{E}$ defined in (\ref{eqn:E}) satisfies
$\mathcal{E} \subset \A$; then, since $\Psi \subset \mathcal{E}$,
$\Psi \setminus \A =\emptyset$. 
Next, since $C \cup D = \X$ and $D$ is closed, it follows that
$\overline{\X \setminus C} \subset D$.  
Then, since $G_{c}(q,z) \subset Q$ for all $(q,z) \in D$, it follows that
for all $(q,z) \in D$,
\begin{equation}
\begin{array}{ccl}
\mu_{V}(q,z) & = & 
\displaystyle V(q,z) - \min_{s \in Q} V(s,z)  \\
& \geq & \displaystyle   V(q,z) -
 \max_{s \in G_{c}(q,z)} V(s,z) \\
& \geq & (1-e^{-1}) V(q,z)  .
\end{array}
\end{equation}
Since $V$ is continuous and positive definite with respect to $\A$,
it follows that $(V,\alpha)$ is an SLFF pair relative to $(\A,\Y)$
with synergy gap exceeding $\varepsilon_{1} V(q,z)$ for any
$\varepsilon_{1} \in (0,1-e^{-1})$.
When $D \cap \A  = \emptyset$, since $V$ is positive definite with
respect to $\A$ and radially unbounded, there exists $\rho>0$ such that
$(q,z) \in D$ implies $V(q,z) \geq \rho$.   In this case, the synergy gap
exceeds any continuous function $\delta$ satisfying
$\delta(q,z) < (1-e^{-1}) 0.5 \left[ \rho + V(q,z) \right]$.  In particular, the
synergy gap exceeds the function $\delta$ given as
$\delta(q,z)= \varepsilon_{1} V(q,z) + \varepsilon_{2}$ where
$\varepsilon_{1} \in (0,0.5(1-e^{-1}))$ and $\varepsilon_{2} \in (0,
0.5 (1-e^{-1}) \rho)$.
\null \hfill \null $\blacksquare$

\subsection{Proof of Theorem \ref{theorem:1}}

Consider the synergistic Lyapunov function
$V$ and feedback $\kappa$.  We claim that
\begin{equation}
(C \setminus \A) \cap [ \Psi \cup \overline{\X \setminus \Y}] = \emptyset .
\label{eq:34}
\end{equation}
Indeed
\begin{equation}
  \mu_{V}(q,z) \leq \delta(q,z)   \qquad \forall (q,z) \in C
\end{equation}
while 
\begin{equation}
\mu_{V}(q,z) > \delta(q,z)  \qquad \forall 
(q,z) \in [\Psi \cup \overline{\X \setminus \Y}] \setminus \A .
\end{equation}
These bounds establish (\ref{eq:34}).

The condition (\ref{eq:34}) together with the fact that $\A \subset \Y \subset \X$ implies
that $C \subset \Y$.  By assumption,
(\ref{eq:derivative}) holds for all $(q,z) \in \Y$
and thus (\ref{eq:derivative}) holds for all $(q,z) \in C$.

By the construction of $D$ and $G_{c}$ in (\ref{eq:controller_b}), 
for all $(q,z) \in D$ and $g_{c} \in G_{c}(z)$, we have 
\begin{equation}
  \begin{array}{rcl}
      V(g_{c},z)  =  \min_{s \in Q} V(s,z)
                       & = & V(q,z) - \mu_{V}(q,z) \\
                       & \leq & V(q,z) - \delta(q,z) .
  \end{array}
\end{equation} 
In particular $V(g_{c},z) - V(q,z) \leq 0$ for all $(q,z) \in D$ and
$g_{c} \in G_{c}(z)$, and $V(g_{c},z) - V(q,z) = 0$ implies $(q,z) \in
\A$.  Using the properties of $V$ and $\delta$, it follows that the
set $\mathcal{A}$ is stable and all solutions are bounded.  It remains
to establish that all complete solutions converge to $\mathcal{A}$.
Note that $\A \subset C$ since $(q,z) \in \A$ implies $\mu_{V}(q,z) =
0 \leq \delta(q,z)$.  Then, by the invariance principle in
\cite{Sanfelice2007}, all complete solutions to (\ref{eq:closed_loop})
converge to the largest weakly invariant set of
\begin{equation}
\left.
\begin{array}{ccl}
    \dot{q} & = & 0 \\
    \dot{z} & = & f(q,z,\kappa(q,z))
\end{array}
\right\}  \quad (q,z) \in \mathcal{E} \cap C  .
\end{equation}
According to the definition of $\Psi$, this weakly invariant set must be contained in $\Psi  \cap C$.
It follows from (\ref{eq:34}) that $\Psi \cap C \subset \A$.
 Thus all complete solutions must converge to $\mathcal{A}$.
\null \hfill \null $\blacksquare$

\subsection{Proof of Theorem \ref{theorem:2blue}}

For all $(q,\zeta) \in \X_1$,
 \begin{equation}
 \renewcommand{\arraystretch}{1.2}
 \begin{array}{l}
\langle \nabla V_1(q,\zeta), \phi_{1}(q,\zeta) + \psi_{1}(q,\zeta)
\kappa_{1}(q,\zeta)  \rangle   \\
 \leq  \langle \nabla V_0(q,z), \phi_{0}(q,z) + \psi_{0}(q,z) \omega \rangle \\
 \quad   - \theta( | \omega-\kappa_{0}(q,z) | )
 - \langle \nabla V_0(q,z),\psi_{0}(q,z)(\omega-
\kappa_{0}(q,z)) \rangle \\
 = \langle \nabla V_0(q,z), \phi_{0}(q,z) + \psi_{0}(q,z) 
\kappa_{0}(q,z) \rangle \\
\quad - \theta( | \omega-\kappa_{0}(q,z) |)
 \leq 0 .
 \end{array}
 \label{eq:hereyblue2}
 \end{equation}
 Define
 \begin{equation}
\renewcommand{\arraystretch}{1.4}
  \begin{array}{rl}
   \mathcal{E}_{1} : = & \left\{ (q,\zeta) \in \X_1 : \right. \\ 
    &\hspace{-0.1in} \left. \langle \nabla V_1(q,\zeta), \phi_{1}(q,\zeta) + \psi_{1}(q,\zeta) \kappa_{1}(q,\zeta) \rangle = 0 \right\} , \\
   \mathcal{W}_{1} : = & \left\{ (q,\zeta) \in \X_{1} :
\langle \nabla V_1(q,\zeta), \psi_{1}(q,\zeta)\rangle = 0 \right\} .
   \end{array}
\end{equation}
Let $\mathcal{E}_{0}$, $\mathcal{W}_{0}$, and $\Omega_{0}$ come from the definitions in
Section \ref{section:SLFF2} for the weak SLFF pair $(V_{0},\kappa_{0})$
for the system (\ref{eq:1aaaa_blue}).
It follows from (\ref{eq:hereyblue2}), the properties of $\theta$, 
the definition of $\psi_{1}$ in (\ref{eq:16blue}), and the
definition of $V_{1}$ in (\ref{eq:V1blue}) that
\begin{equation}
\mathcal{E}_{1} = \left\{
  (q,z) \in \mathcal{E}_{0}  , \ 
  \omega = \kappa_{0}(q,z)   \right\} \subset \mathcal{W}_{1} .
  \label{eq:E1blue}
\end{equation}
Let $\Psi_{1} \subset \X_1$  denote the largest weakly invariant set for the system
\begin{equation}
\left.
  \begin{array}{ccl}
  \dot{q} & = & 0 \\
 \dot{\zeta} & = & \phi_{1}(q,\zeta) + \psi_{1}(q,\zeta) \kappa_{1}(q,\zeta)
   \end{array}
   \right\}  \qquad (q,\zeta) \in \mathcal{E}_{1} .
\end{equation}
It follows from the definition of $\kappa_{1}$ in (\ref{eq:kappa1blue}), the fact that $\dot{\omega}=\kappa_{1}(q,\zeta)$
and the characterization of $\mathcal{E}_{1}$ in (\ref{eq:E1blue}) that
 \begin{equation}
\Psi_{1}  =  \left\{(q,\zeta) \in \X_1:  
(q,z) \in \Omega_{0}  , \ \omega = \kappa_{0}(q,z) \right\} .
\label{eq:Psi1_1}
 \end{equation}
Then, it follows from (\ref{eq:V1blue}) that
\begin{displaymath}
\begin{aligned}
 \mu_{V_{1}}(q,\zeta) &\geq  \mu_{V_{0}}(q,z) \\
&\quad +
 \sigma(\omega-\kappa_{0}(q,z)) - \max_{s \in Q}\sigma(\omega-\kappa_{0}(s,z)).
\end{aligned}
 \end{displaymath}
 Note that $\X_{1} \setminus \Y_{1}= \emptyset$
 and $(q,\zeta) \in \Psi_{1} \setminus \A_{1}$ implies that $(q,z)
 \in \Omega_0 \setminus \A_{0}$.
 Therefore, for $(q,\zeta) \in 
 \left( \Psi_{1}  \cup \overline{\X_{1} \setminus \Y_{1}} \right) \setminus \A_{1}$,
 \begin{displaymath}
  \begin{array}{ccl}
  \mu_{V_{1}}(q,\zeta)
  & \geq & \mu_{V_{0}}(q,z) -
  \max_{s \in Q} \sigma(\kappa_{0}(q,z) - \kappa_{0}(s,z)) \\
  & \geq & \mu_{V_{0}}(q,z) - \varrho(q,z)  \\
  & > & \delta(q,z) .
  \end{array}
\end{displaymath}
Thus, $(V_{1},\kappa_{1})$ is an SLFF pair with gap exceeding
$\delta$.  \null \hfill \null $\blacksquare$
 
 \subsection{Proof of Theorem \ref{theorem:2red}}
 
For all $(q,\zeta) \in \Y_1$,
 \begin{equation}
 \renewcommand{\arraystretch}{1.2}
 \begin{array}{l}
\langle \nabla V_1(q,\zeta), \phi_{1}(q,\zeta) + \psi_{1}(q,\zeta)
\kappa_{1}(q,\zeta)  \rangle \\
 \leq   \langle \nabla V_0(q,z), \phi_{0}(q,z) + \psi_{0}(q,z) \omega \rangle
\\ 
\quad    - \theta( | \omega-\kappa_{0}(q,z) | ) \\
\quad - \langle \nabla V_0(q,z),\psi_{0}(q,z)(\omega-
\kappa_{0}(q,z)) \rangle \\
 = \langle \nabla V_0(q,z), \phi_{0}(q,z) + \psi_{0}(q,z) 
\kappa_{0}(q,z) \rangle \\
\quad - \theta( | \omega-\kappa_{0}(q,z) |)
 \leq 0 .
 \end{array}
 \label{eq:hereyblue}
 \end{equation}
 Define
 \begin{equation}
\renewcommand{\arraystretch}{1.4}
  \begin{array}{rl}
   \mathcal{E}_{1} : = & \left\{ (q,\zeta) \in \Y_1 : \right. \\ 
    &\hspace{-0.1in} \left. \langle \nabla V_1(q,\zeta), \phi_{1}(q,\zeta) + \psi_{1}(q,\zeta) \kappa_{1}(q,\zeta) \rangle = 0 \right\} , \\
   \mathcal{W}_{1} : = & \left\{ (q,\zeta) \in \Y_{1} :
\langle \nabla V_1(q,\zeta), \psi_{1}(q,\zeta)\rangle = 0 \right\} .
   \end{array}
\end{equation}
Let $\mathcal{E}_{0}$, $\mathcal{W}_{0}$, and $\Omega_{0}$ come from the definitions in
Section \ref{section:SLFF2} for the weak SLFF pair $(V_{0},\kappa_{0})$
for the system (\ref{eq:1aaaa_blue}).
It follows from (\ref{eq:hereyblue}), the properties of $\theta$, 
the definition of $\psi_{1}$ in (\ref{eq:16blue}), and the
definition of $V_{1}$ in (\ref{eq:V1red}) that
\begin{equation}
\mathcal{E}_{1} = \left\{
  (q,z) \in \mathcal{E}_{0}  , \ 
  \omega = \kappa_{0}(q,z)   \right\} \subset \mathcal{W}_{1} .
  \label{eq:E1}
\end{equation}
Let $\Psi_{1} \subset \X_1$  denote the largest weakly invariant set for the system
\begin{equation}
\left.
  \begin{array}{ccl}
  \dot{q} & = & 0 \\
 \dot{\zeta} & = & \phi_{1}(q,z) + \psi_{1}(q,z) \kappa_{1}(q,z)
   \end{array}
   \right\}  \qquad (q,\zeta) \in \mathcal{E}_{1} .
\end{equation}
It follows from the definition of $\kappa_1$ in (\ref{eq:kappa1red}), the fact that $\dot{\omega}=\kappa_{1}(q,\zeta)$
and the characterization of $\mathcal{E}_{1}$ in (\ref{eq:E1}) that
 \begin{equation}
\Psi_{1}  =  \left\{(q,\zeta) \in \X_1:  
(q,z) \in \Omega_{0}  , \ \omega = \kappa_{0}(q,z) \right\} .
\label{eq:Psi1_2}
 \end{equation}
 Note that $(q,\zeta) \in \Psi_{1} \setminus \A_{1}$ implies
 $(q,z) \in \Omega_{0} \setminus \A_{0}$.
 Also $\overline{\X_{1} \setminus \Y_{1}}
 =\overline{\X_{0} \setminus \Y_{0}} \times \Re^{m}$ so
 that 
 \begin{displaymath}
 \left( \Psi_{1} \cup \overline{ \X_{1} \setminus \Y_{1}} \right)
\setminus \A_{1} 
\subset \left( (\Omega_{0} \setminus \A_{0}) \cup \overline{\X_{0} \setminus
\Y_{0}} \right) \times \Re^{m} .
 \end{displaymath}
 Then, it follows from (\ref{eq:V1red}) and the facts that $\rho \circ \sigma_{2}$
 is globally Lipschitz with Lipschitz constant less than or equal to $L>0$ and
 $V_{0}$ is type I
 ready-made relative to $\sigma$  with $\sigma(v):=L|v|$ for all $v \in \Re^{m}$
 that,
for $(q,\zeta) \in \left( \Psi_{1} \cup \overline{ \X_{1} \setminus \Y_{1}} \right)
\setminus \A_{1}$
 \begin{displaymath}
 \renewcommand{\arraystretch}{1.4}
 \begin{array}{rcl}
 \mu_{V_{1}}(q,\zeta) & \geq &  \displaystyle \mu_{V_{0}}(q,z) \\
& & \displaystyle  \! \! \!  + \rho(\sigma_{2}(\omega-\kappa_{0}(q,z))) - \max_{s \in Q}
\rho(\sigma_{2}(\omega-\kappa_{0}(s,z))) \\
& \geq & \displaystyle \mu_{V_{0}}(q,z) - \max_{s \in Q}
L |\kappa_{0}(s,z) - \kappa_{0}(q,z)| \\
  & \geq & \mu_{V_{0}}(q,z) - \varrho(q,z)  \\
  & > & \delta(q,z) .
\end{array}
 \end{displaymath}
Thus, $(V_{1},\kappa_{1})$ is an SLFF pair with gap exceeding
$\delta$.  \null \hfill \null $\blacksquare$

\end{document}